\documentclass[12pt]{article} 
\def\Bbb{\bf} 
\newcommand\C{{ \Bbb C}}

\newcommand\PP{{\Bbb P}}

\newtheorem{thm}{Theorem}

\def\comment#1{ }

\begin{document}
\title{A $W(E_6)$-equivariant projective embedding of the moduli space of cubic surfaces}
\author{Masaaki YOSHIDA \thanks{Department of Mathematics, Kyushu University, Fukuoka 812 Japan}}
\date{22 December 1999}
\maketitle
\begin{abstract}
An explicit projective embedding of the moduli space of marked cubic surfaces is given. This embedding is equivariant under the Weyl group of type $E_6$. The image is defined by a system of linear and cubic equations. To express the embedding in a most symmetric way, the target would be 79-dimensional, however the image lies in a 9-dimensional linear subspace.

\end{abstract}
\par\smallskip
\section{The moduli space of cubic surfaces}
\subsection{The space $M$ and the action of the group $\underline G$}
We first fix some notation and recall a few known facts on the moduli space
of marked cubic surfaces.
The moduli space of marked cubic surfaces, which we denote by $M$,
is studied for example in \cite{Na} and \cite{SeYo}. 
Since any nonsingular cubic
surface can be obtained by blowing up the projective plane $\PP^2$ at
six points, it can be represented by a $3\times6$-matrix of which
columns give homogeneous coordinates of the six points. In order to
get a smooth cubic surface from six points, 
we assume that no three points are
collinear and the six points are not on a conic. 
On the set of $3\times 6$ matrices, we have a cannical action
of $GL_3$ on the left and the group $\C^\times$ acts
naturally on homogeneous coordinates. 
By killing such ambiguity of coordinates, we get the following expression
\[ x=\left( 
\begin{array}{cccccc}
1&0&0&1& 1& 1\\ 0&1&0&1&x_1&x_2\\ 0&0&1&1&x_3&x_4
\end{array}
\right); \] 
in this paper we use local coordinates $(x_1,x_2,x_3,x_4)$ on $M$.
The six points represented by the matrix above produces a non-singular cubic surface if and only if the following quantitiy does not vanish.
\begin{eqnarray*} D(x)&:=
&x_1x_2x_3x_4(x_1-1)(x_2-1)(x_3-1)(x_4-1)\\
&&\times(x_1-x_2)(x_1-x_3)(x_2-x_4)(x_3-x_4)D_1D_2Q,
\end{eqnarray*}
where
\begin{eqnarray*}
D_1&:=&x_1x_4-x_2x_3,\\ 
D_2&:=&x_1x_4-x_4+x_2-x_2x_3+x_3-x_1,\\
Q&:=&-x_2x_3x_1-x_2x_3x_4+x_2x_3+x_1x_4x_2+x_1x_4x_3-x_1x_4,\\
\end{eqnarray*}
Thus we can identify the moduli space $M$ 
with the affine open set $\{x=(x_1,\dots,x_4)\mid D(x)\not=0\}.$

Let us define as in \cite{SeYo} six bi-rational transformations $s_1,\dots,s_6$ in 
$x=(x_1,\dots,x_4):$
\begin{eqnarray*}
s_1:(x_1,x_2,x_3,x_4) & \rightarrow &
 \left({1\over x_1},{1\over x_2},{x_3\over x_1},{x_4\over x_2}\right), \\
s_2:(x_1,x_2,x_3,x_4) & \rightarrow &
  (x_3,x_4,x_1,x_2) ,\\
s_3:(x_1,x_2,x_3,x_4) & \rightarrow &
  \left({x_1-x_3\over 1-x_3}, {x_2-x_4\over 1-x_4},
        {x_3\over x_3-1}, {x_4\over x_4-1}\right) ,\\
s_4:(x_1,x_2,x_3,x_4) & \rightarrow &
  \left({1\over x_1},{x_2\over x_1}, {1\over x_3}, {x_4\over x_3}\right),\\
s_5:(x_1,x_2,x_3,x_4) & \rightarrow &
  (x_2, x_1, x_4, x_3) ,\\
s_6:(x_1,x_2,x_3,x_4) & \rightarrow &
  \left({1\over x_1}, {1\over x_2}, {1\over x_3}, {1\over x_4}\right).
\end{eqnarray*}
If $M$ is regarded as the configuration space of six points in $\PP^2$, 
the transformation $s_1$, for example, corresponds to the interchange of 
the two points represented by the first two column vectors of the matrix $x$.
Each $s_i$ turns out to be a bi-regular involution on $M$, 
and they form a group $\underline G$ 
isomorphic to the Weyl group of type $E_6$; relation of the generators are
given by the Coxeter graph
$$
\begin{array}{ccccccccc}
s_1&\mbox{---}&s_2&\mbox{---}&s_3&\mbox{---}&s_4&\mbox{---}&s_5\\
                   &  &   &  &\mid& &   &  &   \\
                   &  &   &  &s_6&  &   &  &        
\end{array}$$
\subsection{Root system $\Delta$ of type $E_6$}
We review the root system $\Delta$ of type $E_6$, following \cite{H}.
Consider an 8-dimensional Euclidean space 
${\tilde E}$ with
a standard basis $\varepsilon_1,\cdots,\varepsilon_8$.
Let $\langle\cdot,\cdot\rangle$ be the inner product on ${\tilde E}$
defined by
$\langle\varepsilon_j, \varepsilon_k\rangle = \delta_{jk}$
and let $E$ be the linear subspace of ${\tilde E}$
spanned by the six vectors 
\[\varepsilon_1,\cdots,\varepsilon_5,\>
{\tilde \varepsilon}=\varepsilon_6-\varepsilon_7-\varepsilon_8.\]
We introduce the 36 vectors:
\[
r=-{1 \over 2}(
\varepsilon_1+\varepsilon_2+\varepsilon_3+\varepsilon_4+\varepsilon_5
+{\tilde \varepsilon}),\]
\[r_{1j}=-\varepsilon_{j-1}+r_0, \quad 2 \le   j\le 6\]
\[r_{jk}= \varepsilon_{j-1}-\varepsilon_{k-1}, \quad 2 \le   j<  k\le 6\]
\[r_{1jk}= -\varepsilon_{j-1}-\varepsilon_{k-1}, \quad 2 \le   j<  k\le 6\]
\[r_{ijk}= -\varepsilon_{i-1}-\varepsilon_{j-1}-\varepsilon_{k-1}+r_0, 
\quad 2 \le  i<j<k\le 6\]
where $r_{ij}=r_{ji}, r_{ijk}=r_{jik}=r_{ikj},$
\[r_0={1 \over 2}(\varepsilon_1+\varepsilon_2+\varepsilon_3
+\varepsilon_4+\varepsilon_5-{\tilde \varepsilon}).\]
Note that 
$$r\perp r_{ij},\quad r_{ij}\perp r_{kl},\quad r_{ij}\perp r_{ijk},\quad r_{ij}\perp r_{klm},\quad r_{ijk}\perp r_{ilm}.$$
The set 
$$\Delta=\{\pm r,\ \pm r_{ij},\ \pm r_{ijk}\}$$
forms a root system of type $E_6$.
For example,
$$r_{12},\quad r_{123},\quad r_{23},\quad r_{34},\quad r_{45},\quad r_{56}$$
can serve as a system  of positive simple roots; its extended Dynkin diagram is given as

\begin{center}\begin{tabular}
{p{5mm} p{10mm} p{5mm} p{10mm} p{5mm} p{10mm} p{5mm} p{10mm} p{5mm}}
$r_{12}$ &------&$r_{23}$&------&$r_{34}$&------&$r_{45}$
&------&$r_{56}$\\
      &      &       &      & $\mid$&      &       &      & \\
      &      &       &      & $r_{123}$ &  &       &      & \\
      &      &       &      & $\mid$&      &       &      & \\
      &      &       &      & $r$ &  &       &      & \\
\end{tabular}\end{center}
The set 
$\{r$, $r_{jk}$, $r_{ijk}\}$ is the totality of
positive roots of $\Delta$.

Let $s_r, s_{ij}$ and $s_{ijk}$ be the reflections on $E$ with respect to $r, r_{ij}$ and $r_{ijk}$. These reflections act on $\Delta$ as 
$$\begin{array}{ll}\quad\ \  s_r: &r_{ij}\leftrightarrow r_{ij},\quad r_{ijk}\leftrightarrow r_{lmn},\quad 
\{i,j,k,l,m,n\}=\{1,\dots,6\},\\
\quad\ \  s_{ij}: &{\rm permutation\ of\ the\ indices\ }i {\rm \ and\ }j,\end{array}$$
$$\begin{array}{llllllllll}
s_{123}:&r_{12}&\leftrightarrow &r_{12},\quad &r_{14}&\leftrightarrow &r_{234},\quad &r_{56}&\leftrightarrow &r_{56},\\
&r_{145}&\leftrightarrow &r_{145},\quad &r&\leftrightarrow &r_{456},&&&\end{array}$$
modulo signs. Let us define two reflection group
$$
G_1=\langle s_{12},s_{23},s_{34},s_{45},s_{56}\rangle\cong S_6,\quad
G=\langle G_1,s_{123}\rangle\cong W(E_6),$$
where $S_6$ is the symmetric group on six numerals $\{1,\dots,6\}$, $W(E_6)$ the Weyl group of type $E_6$, and $\langle a,b,\dots\rangle$ denotes the group generated by $a,b,\dots$   Note that
$$G_1\subset \langle G_1,s_r\rangle= G_1\times\langle s_r\rangle
\subset G,$$
and $G$ acts transitively on $\Delta$.
\subsection{Naruki's cross-ratio variety}
A smooth compactification of $M$ known as Naruki's cross-ratio
variety $\mathcal C$ (\cite{Na}, \cite{SeYo}), embedded in $(\PP^1)^{45}$, is the union of $M$ and the 76 divisors. The 36 of them correspond 
to the positive roots of $\Delta$ (they are said to be of the first kind), the other 40 divisors (said to be of the second kind, and are isomorphic to $(\PP^1)^3$) can be represented
as follows:
Take three subsets
$\Delta_1$, $\Delta_2$, $\Delta_3$ of $\Delta$ satisfying the following conditions:
\begin{itemize}
\item Each of 
$\Delta_1$, $\Delta_2$, $\Delta_3$ is a root system of type $A_2$.
\item $\Delta_1$, $\Delta_2$, $\Delta_3$ are mutually orthogonal.
\item The vectors in
$\Delta_1\cup\Delta_2 \cup \Delta_3$ span $E$.
\end{itemize}
Note that each one of such three root systems determines the other two.

 Such a triple  $\{\Delta_1,\Delta_2,\Delta_3\}$ dertermines a divisor. 
According to the naming of the roots, we must use two different expressions:
The first one is of the form
\[
\{\pm r_{12},\pm r_{23},\pm r_{13}\},
\quad
\{\pm r_{45},\pm r_{56},\pm r_{46}\},
\quad
\{\pm r,\pm r_{123},\pm r_{456}\},
\]
(the corresponding divisor is denoted by $Z_{123,456}$ in \cite{SeYo},) and
the second one is of the form
\[
\{\pm r_{12},\pm r_{234},\pm r_{134}\},
\quad
\{\pm r_{34},\pm r_{356},\pm r_{456}\},
\quad
\{\pm r_{56},\pm r_{125},\pm r_{126}\},
\]
(the corresponding divisor is denoted by $Z_{12,34,56}$ in \cite{SeYo}).
Note that
 $$Z_{123,456}=Z_{456,123},\quad
Z_{12,34,56}=Z_{56,12,34}\not=Z_{12,56,34}.
$$
Thus, permuting the indices under $G_1$, we have 40 ($=10+30$) such divisors.
The group $G$ acts transitively on these 40 divisors.
\par\medskip\noindent
{\bf Remark:}
These divisors (of the second kind) are disjoint to each other, and can be blown-down to points. In fact they correspond bijectively to the cusps of the modular group studied in
 \cite{ACT} (see also \cite{Ma2} and \cite{SaYo2}). 
\section{Embedding $\varphi:M\to \PP^{80-1}$}
\subsection{Coordinates on $\PP^{80-1}$ and the action of $G$}
Let $\mathcal A$ be the set of 40 labels $(123,456)$ and $(12,34,56)$
with the following identification
\begin{eqnarray*}
(123,456)&=&(213,456)=(132,456)=(456,123),\\
(12,34,56)&=&(21,34,56)=(56,12,34).
\end{eqnarray*}
Since $G$ acts on the set $\Delta$ of roots, it acts also on the set of
40 divisors above, and so that it also acts on the set $\mathcal A$ of
40 labels.

We introduce 80 homogeneous coordinates 
$$y_\alpha,\quad y_{-\alpha},\qquad \alpha\in\mathcal A$$
on $\PP^{80-1}$. I define an action of $G$ on $\PP^{80-1}$ by the
following action of the generators $s_{12},\dots,s_{56}$ and $s_{123}$
on the coordinates. Let $s$ be one of the generators and $\alpha\in
\mathcal A$; we assign
\begin{eqnarray*}
s(y_\alpha)&=&y_\alpha,\quad s(y_{-\alpha})=y_{-\alpha}\qquad \mbox{if
 $s\alpha=\alpha$,}\\
s(y_\alpha)&=&y_{-\beta},\quad s(y_{-\alpha})=y_{\beta}\qquad \mbox{if
 $s\alpha=\beta\not=\alpha$.}
\end{eqnarray*}

\subsection{Definition of $\varphi$}
In this section we define a map $M\to \PP^{80-1}$. For a $3\times6$
matrix $x=(x_{ij})$, we consider 80 polynomials  of degree 18 as follows:
\begin{eqnarray*}
y_{(123,456)}(x)&=& D_{123}(x)D_{456}(x)Q(x),\\
y_{(12,34,56)}(x)&=& D_{134}(x)D_{234}(x)D_{356}(x)D_{456}(x)D_{512}(x)D_{612}(x),
\end{eqnarray*}
and $y_{-\alpha}(x)=-y_\alpha\ (\alpha\in\mathcal A)$, where $Q(x)$ is the determinant of the $6\times6$-matrix with columns 
$$(x_{1j}x_{2j},x_{2j}x_{3j},x_{3j}x_{1j}, x_{1j}^2,x_{2j}^2,x_{3j}^2)\quad j=1,\dots,6.$$
Since we have
$$y_\alpha(gxh)=(\det g)^6y_\alpha(x)(\det h)^3,$$
the correswpondence above defines a map $\varphi:M\to \PP^{80-1}$.
For later use, we present 40 polynomials $y_\alpha(x)$ in terms of the
coordinates  $(x_1,x_2,x_3,x_4)$ introduced in \S1; the remaining 40
polynomoals are given by $y_{-\alpha}(x)=-y_\alpha(x)$. In the following 
table, $y_\alpha(x)$ is denoted simply by $\alpha$, and we number them as $y_1,\dots,y_{40}$:
\begin{eqnarray*}
y_{1}=(156,234)&:=&D_1Q:\quad y_{2}=(123,456):=D_2Q:\\
y_{3}=(124,356)&:=&(x_2-x_1)Q:\quad y_{4}=(145,236):=(x_3-x_1)Q:\\
y_{5}=(146,235)&:=&(x_4-x_2)Q:\quad y_{6}=(134,256):=(x_4-x_3)Q:\\
y_{7}=(135,246)&:=&x_1(x_4-1)Q:\quad y_{8}=(136,245):=x_2(x_3-1)Q:\\
y_{9}=(125,346)&:=&x_3(x_2-1)Q:\quad y_{10}=(126,345):=x_4(x_1-1)Q:\\
y_{11}=(12,56,34)&:=&D_1(x_1-1)(x_2-1)(x_4-x_3):\\
y_{12}=(16,23,45)&:=&D_1(x_1-1)(x_3-1)(x_4-x_2):\\
y_{13}=(15,23,46)&:=&D_1(x_2-1)(x_4-1)(x_3-x_1):\\
y_{14}=(13,56,24)&:=&D_1(x_3-1)(x_4-1)(x_2-x_1):\\
y_{15}=(15,24,36)&:=&-D_1x_1(x_2-1)(x_3-1):\\
y_{16}=(16,24,35)&:=&-D_1x_2(x_1-1)(x_4-1):\\
y_{17}=(15,34,26)&:=&-D_1x_3(x_1-1)(x_4-1):\\
y_{18}=(16,34,25)&:=&-D_1x_4(x_2-1)(x_3-1):\\
y_{19}=(12,36,45)&:=&D_2x_2x_3(x_1-1):\\
y_{20}=(12,35,46)&:=&D_2x_1x_4(x_2-1):\\
y_{21}=(13,26,45)&:=&D_2x_1x_4(x_3-1):\\
y_{22}=(13,25,46)&:=&D_2x_2x_3(x_4-1):\\
y_{23}=(13,24,56)&:=&D_2x_1x_2(x_4-x_3):\\
y_{24}=(15,46,23)&:=&D_2x_1x_3(x_4-x_2):\\
y_{25}=(16,45,23)&:=&D_2x_2x_4(x_3-x_1):\\
y_{26}=(12,34,56)&:=&D_2x_3x_4(x_2-x_1):\\
y_{27}=(13,46,25)&:=&x_1(x_2-1)(x_3-1)(x_4-x_2)(x_4-x_3):\\
y_{28}=(13,45,26)&:=&x_2(x_1-1)(x_4-1)(x_3-x_1)(x_4-x_3):\\
y_{29}=(12,46,35)&:=&x_3(x_1-1)(x_4-1)(x_2-x_1)(x_4-x_2):\\
y_{30}=(12,45,36)&:=&x_4(x_2-1)(x_3-1)(x_2-x_1)(x_3-x_1):\\
y_{31}=(14,35,26)&:=&-x_1(x_1-1)(x_4-x_2)(x_4-x_3):\\
y_{32}=(14,36,25)&:=&-x_2(x_2-1)(x_3-x_1)(x_4-x_3):\\
y_{33}=(14,25,36)&:=&-x_3(x_3-1)(x_2-x_1)(x_4-x_2):\\
y_{34}=(14,26,35)&:=&-x_4(x_4-1)(x_2-x_1)(x_3-x_1):\\
y_{35}=(16,25,34)&:=&-x_2x_3(x_1-1)(x_4-x_2)(x_4-x_3):\\
y_{36}=(15,26,34)&:=&-x_1x_4(x_2-1)(x_3-x_1)(x_4-x_3):\\
y_{37}=(16,35,24)&:=&-x_1x_4(x_3-1)(x_2-x_1)(x_4-x_2):\\
y_{38}=(15,36,24)&:=&-x_2x_3(x_4-1)(x_2-x_1)(x_3-x_1):\\
y_{39}=(14,56,23)&:=&D_1D_2:\\
y_{40}=(14,23,56)&:=&(x_2-x_1)(x_3-x_1)(x_4-x_2)(x_4-x_3):\\
\end{eqnarray*}

\subsection{$G$-Equivariance of $\varphi$}
Recall that the group $\underline G$ acts on $M$, and that $G$ acts on $\PP^{80-1}$. Let us identify the groups $\underline G$ and $G$ by
$$\iota:s_{12}\mapsto s_1,\dots,s_{56}\mapsto s_5, \ s_{123}\mapsto
s_6.$$ Then we have 
\begin{thm}The map $\varphi:M\to\PP^{80-1}$ is $G$-equivariant:
$$g(\varphi(x))=\varphi(\iota(g)x),\quad g\in G,\ x\in M,$$
that is,
$$(gy_\alpha)(x)=c_gy_\alpha(\iota(g)x),\quad g\in G,\ 
\alpha\in\pm\mathcal A,\ x\in M,$$
where $c_g$ is a rational function in $(x_1,x_2,x_3,x_4)$.\end{thm}
\par\noindent
{\bf Convention:} Once this theorem is established, we ignore the redundant ones $y_{-\alpha}(x)=-y_{\alpha}(x)$ and regard $\varphi$ as the map
$$M\ni x\longmapsto :y_\alpha(x):\in \PP^{40-1}.$$
The group $G$ still acts on $\PP^{40-1}$ by the transformations given in\S2.3.
\par\smallskip
In order to prove the theorem, we have only to check the identity for a set of gnerators of $G$.
Under $s_1$, the fourty polynomials are transformed as follows:
$$\begin{array}{lllll}
y_{1}\to -c_1y_{6},&
y_{2}\to c_1 y_{2},&
y_{3}\to c_1 y_{3},&
y_{4}\to -c_1y_{8},&
y_{5}\to -c_1y_{7},\\
y_{6}\to -c_1y_{1},&
y_{7}\to -c_1y_{5},&
y_{8}\to -c_1y_{4},&
y_{9}\to c_1 y_{9},&
y_{10}\to c_1 y_{10},\\
y_{11}\to c_1 y_{11},&
y_{12}\to -c_1y_{28},&
y_{13}\to -c_1y_{27},&
y_{14}\to -c_1y_{40},&
y_{15}\to -c_1y_{32},\\
y_{16}\to -c_1y_{31},&
y_{17}\to -c_1y_{35},&
y_{18}\to -c_1y_{36},&
y_{19}\to c_1 y_{19},&
y_{20}\to c_1 y_{20},\\
y_{21}\to -c_1y_{25},&
y_{22}\to -c_1y_{24},&
y_{23}\to -c_1y_{39},&
y_{24}\to -c_1y_{22},&
y_{25}\to -c_1y_{21},\\
y_{26}\to c_1 y_{26},&
y_{27}\to -c_1y_{13},&
y_{28}\to -c_1y_{12},&
y_{29}\to c_1 y_{29},&
y_{30}\to c_1 y_{30},\\
y_{31}\to -c_1y_{16},&
y_{32}\to -c_1y_{15},&
y_{33}\to -c_1y_{38},&
y_{34}\to -c_1y_{37},&
y_{35}\to -c_1y_{17},\\
y_{36}\to -c_1y_{18},&
y_{37}\to -c_1y_{34},&
y_{38}\to -c_1y_{33},&
y_{39}\to -c_1y_{23},&
y_{40}\to -c_1y_{14},
\end{array}$$
where $c_1=(x_1x_2)^{-3}$, 
under $s_{2}$,
$$\begin{array}{lllll}
y_{1}\to y_{1},&
y_{2}\to y_{2},&
y_{3}\to -y_{6},&
y_{4}\to y_{4},&
y_{5}\to y_{5},\\
y_{6}\to -y_{3},&
y_{7}\to -y_{9},&
y_{8}\to -y_{10},&
y_{9}\to -y_{7},&
y_{10}\to -y_{8},\\
y_{11}\to -y_{14},&
y_{12}\to y_{12},&
y_{13}\to y_{13},&
y_{14}\to -y_{11},&
y_{15}\to -y_{17},\\
y_{16}\to -y_{18},&
y_{17}\to -y_{15},&
y_{18}\to -y_{16},&
y_{19}\to -y_{21},&
y_{20}\to -y_{22},\\
y_{21}\to -y_{19},&
y_{22}\to -y_{20},&
y_{23}\to -y_{26},&
y_{24}\to y_{24},&
y_{25}\to y_{25},\\
y_{26}\to -y_{23},&
y_{27}\to -y_{29},&
y_{28}\to -y_{30},&
y_{29}\to -y_{27},&
y_{30}\to -y_{28},\\
y_{31}\to -y_{33},&
y_{32}\to -y_{34},&
y_{33}\to -y_{31},&
y_{34}\to -y_{32},&
y_{35}\to -y_{37},\\
y_{36}\to -y_{38},&
y_{37}\to -y_{35},&
y_{38}\to -y_{36},&
y_{39}\to y_{39},&
y_{40}\to y_{40},\end{array}$$
under $s_3$,
$$\begin{array}{lllll}
y_{1}\to c_3y_{1},&
y_{2}\to -c_3y_{3},&
y_{3}\to -c_3y_{2},&
y_{4}\to -c_3y_{7},&
y_{5}\to -c_3y_{8},\\
y_{6}\to c_3y_{6},&
y_{7}\to -c_3y_{4},&
y_{8}\to -c_3y_{5},&
y_{9}\to c_3y_{9},&
y_{10}\to c_3y_{10},\\
y_{11}\to c_3y_{11},&
y_{12}\to -c_3y_{16},&
y_{13}\to -c_3y_{15},&
y_{14}\to -c_3y_{39},&
y_{15}\to -c_3y_{13},\\
y_{16}\to -c_3y_{12},&
y_{17}\to c_3y_{17},&
y_{18}\to c_3y_{18},&
y_{19}\to -c_3y_{29},&
y_{20}\to -c_3y_{30},\\
y_{21}\to -c_3y_{34},&
y_{22}\to -c_3y_{33},&
y_{23}\to -c_3y_{40},&
y_{24}\to -c_3y_{38},&
y_{25}\to -c_3y_{37},\\
y_{26}\to c_3y_{26},&
y_{27}\to -c_3y_{32},&
y_{28}\to -c_3y_{31},&
y_{29}\to -c_3y_{19},&
y_{30}\to -c_3y_{20},\\
y_{31}\to -c_3y_{28},&
y_{32}\to -c_3y_{27},&
y_{33}\to -c_3y_{22},&
y_{34}\to -c_3y_{21},&
y_{35}\to c_3y_{35},\\
y_{36}\to c_3y_{36},&
y_{37}\to -c_3y_{25},&
y_{38}\to -c_3y_{24},&
y_{39}\to -c_3y_{14},&
y_{40}\to -c_3y_{23},\end{array}$$
where $c_3=(1-x_3)^{-3}(1-x_4)^{-3}$, under $s_4$,
$$\begin{array}{lllll}
y_{1}\to -c_4y_{5},&
y_{2}\to c_4y_{2},&
y_{3}\to -c_4y_{9},&
y_{4}\to c_4y_{4},&
y_{5}\to -c_4y_{1},\\
y_{6}\to -c_4y_{7},&
y_{7}\to -c_4y_{6},&
y_{8}\to c_4y_{8},&
y_{9}\to -c_4y_{3},&
y_{10}\to c_4y_{10},\\
y_{11}\to -c_4y_{29},&
y_{12}\to c_4y_{12},&
y_{13}\to -c_4y_{40},&
y_{14}\to -c_4y_{27},&
y_{15}\to -c_4y_{33},\\
y_{16}\to -c_4y_{35},&
y_{17}\to -c_4y_{31},&
y_{18}\to -c_4y_{37},&
y_{19}\to c_4y_{19},&
y_{20}\to -c_4y_{26},\\
y_{21}\to c_4y_{21},&
y_{22}\to -c_4y_{23},&
y_{23}\to -c_4y_{22},&
y_{24}\to -c_4y_{39},&
y_{25}\to c_4y_{25},\\
y_{26}\to -c_4y_{20},&
y_{27}\to -c_4y_{14},&
y_{28}\to c_4y_{28},&
y_{29}\to -c_4y_{11},&
y_{30}\to c_4y_{30},\\
y_{31}\to -c_4y_{17},&
y_{32}\to -c_4y_{38},&
y_{33}\to -c_4y_{15},&
y_{34}\to -c_4y_{36},&
y_{35}\to -c_4y_{16},\\
y_{36}\to -c_4y_{34},&
y_{37}\to -c_4y_{18},&
y_{38}\to -c_4y_{32},&
y_{39}\to -c_4y_{24},&
y_{40}\to -c_4y_{13},\end{array}$$
where $c_4=(x_1x_3)^{-3}$, under $s_5$,
$$\begin{array}{lllll}
y_{1}\to y_{1},&
y_{2}\to y_{2},&
y_{3}\to y_{3},&
y_{4}\to -y_{5},&
y_{5}\to -y_{4},\\
y_{6}\to y_{6},&
y_{7}\to -y_{8},&
y_{8}\to -y_{7},&
y_{9}\to -y_{10},&
y_{10}\to -y_{9},\\
y_{11}\to y_{11},&
y_{12}\to -y_{13},&
y_{13}\to -y_{12},&
y_{14}\to y_{14},&
y_{15}\to -y_{16},\\
y_{16}\to -y_{15},&
y_{17}\to -y_{18},&
y_{18}\to -y_{17},&
y_{19}\to -y_{20},&
y_{20}\to -y_{19},\\
y_{21}\to -y_{22},&
y_{22}\to -y_{21},&
y_{23}\to y_{23},&
y_{24}\to -y_{25},&
y_{25}\to -y_{24},\\
y_{26}\to y_{26},&
y_{27}\to -y_{28},&
y_{28}\to -y_{27},&
y_{29}\to -y_{30},&
y_{30}\to -y_{29},\\
y_{31}\to -y_{32},&
y_{32}\to -y_{31},&
y_{33}\to -y_{34},&
y_{34}\to -y_{33},&
y_{35}\to -y_{36},\\
y_{36}\to -y_{35},&
y_{37}\to -y_{38},&
y_{38}\to -y_{37},&
y_{39}\to y_{39},&
y_{40}\to y_{40},\end{array}$$
and under $s_6$,
$$\begin{array}{lllll}
y_{1}\to-c_6y_{39},&
y_{2}\to c_6y_{2},&
y_{3}\to-c_6y_{26},&
y_{4}\to-c_6y_{25},&
y_{5}\to-c_6y_{24},\\
y_{6}\to-c_6y_{23},&
y_{7}\to-c_6y_{22},&
y_{8}\to-c_6y_{21},&
y_{9}\to-c_6y_{20},&
y_{10}\to-c_6y_{19},\\
y_{11}\to c_6y_{11},&
y_{12}\to c_6y_{12},&
y_{13}\to c_6y_{13},&
y_{14}\to c_6y_{14},&
y_{15}\to-c_6y_{18},\\
y_{16}\to-c_6y_{17},&
y_{17}\to-c_6y_{16},&
y_{18}\to-c_6y_{15},&
y_{19}\to-c_6y_{10},&
y_{20}\to-c_6y_{9},\\
y_{21}\to-c_6y_{8},&
y_{22}\to-c_6y_{7},&
y_{23}\to-c_6y_{6},&
y_{24}\to-c_6y_{5},&
y_{25}\to-c_6y_{4},\\
y_{26}\to-c_6y_{3},&
y_{27}\to c_6y_{27},&
y_{28}\to c_6y_{28},&
y_{29}\to c_6y_{29},&
y_{30}\to c_6y_{30},\\
y_{31}\to-c_6y_{35},&
y_{32}\to-c_6y_{36},&
y_{33}\to-c_6y_{37},&
y_{34}\to-c_6y_{38},&
y_{35}\to-c_6y_{31},\\
y_{36}\to-c_6y_{32},&
y_{37}\to-c_6y_{33},&
y_{38}\to-c_6y_{34},&
y_{39}\to-c_6y_{1},&
y_{40}\to c_6y_{40},
\end{array}$$
where $c_6=(x_1x_2x_3x_4)^{-2}$.
Maybe it is interesting to see what happens under the operation of the 
involution $s_r$ (classically called the association) which sends $(x_1,x_2,x_3,x_4)$ to
$$\left(
{(x_4-1)D_1\over(x_4-x_2)(x_4-x_3)},
{(x_3-1)D_1\over(x_3-x_1)(x_4-x_3)},
{(x_2-1)D_1\over(x_4-x_2)(x_2-x_1)},
{(x_1-1)D_1\over(x_3-x_1)/(x_2-x_1)}\right):$$
$$\begin{array}{lllll}
y_{1}\to c_ry_{1},&
y_{2}\to c_ry_{2},&
y_{3}\to c_ry_{3},&
y_{4}\to c_ry_{4},&
y_{5}\to c_ry_{5},\\
y_{6}\to c_ry_{6},&
y_{7}\to c_ry_{7},&
y_{8}\to c_ry_{8},&
y_{9}\to c_ry_{9},&
y_{10}\to c_ry_{10},\\
y_{11}\to-c_ry_{26},&
y_{12}\to-c_ry_{25},&
y_{13}\to-c_ry_{24},&
y_{14}\to-c_ry_{23},&
y_{15}\to-c_ry_{38},\\
y_{16}\to-c_ry_{37},&
y_{17}\to-c_ry_{36},&
y_{18}\to-c_ry_{35},&
y_{19}\to-c_ry_{30},&
y_{20}\to-c_ry_{29},\\
y_{21}\to-c_ry_{28},&
y_{22}\to-c_ry_{27},&
y_{23}\to-c_ry_{14},&
y_{24}\to-c_ry_{13},&
y_{25}\to-c_ry_{12},\\
y_{26}\to-c_ry_{11},&
y_{27}\to-c_ry_{22},&
y_{28}\to-c_ry_{21},&
y_{29}\to-c_ry_{20},&
y_{30}\to-c_ry_{19},\\
y_{31}\to-c_ry_{34},&
y_{32}\to-c_ry_{33},&
y_{33}\to-c_ry_{32},&
y_{34}\to-c_ry_{31},&
y_{35}\to-c_ry_{18},\\
y_{36}\to-c_ry_{17},&
y_{37}\to-c_ry_{16},&
y_{38}\to-c_ry_{15},&
y_{39}\to-c_ry_{40},&
y_{40}\to-c_ry_{39},
\end{array}
$$
where 
$$c_r=\left({D_1D_2\over(-x_4+x_2)(x_4-x_3)(x_1-x_3)(x_1-x_2)}\right)^3.$$
\subsection{$\varphi$ embeds $M$}
It is known in \cite{Yo3} that the map
$$M\ni x\longrightarrow y_1(x):y_3(x):y_4(x):y_5(x):y_7(x)\in \PP^4$$
is two-to-one, and induces an embedding of
the quotient space $M/\langle s_r\rangle$.
Thus the composite of $\varphi$ and the projection
$$M\longrightarrow\PP^{40-1}\longrightarrow\PP^4$$
is a two-to-one map.
This fact together with the equivariance of $\varphi$ under the involution $s_r$ shown just above
implies 
\begin{thm} $\varphi$ embeds $M$ into $\PP^{40-1}$. \end{thm}
\subsection{Prolongation of $\varphi$ to degenerate arrangements}
Let us consider degenerate arrangements of six points on the plane. Since arrangements with three collinear points can be transformed under $G$ to those with six points on a conic, we assume, Without loss
of generality, that our arrangements represented by $x=(x_1,\dots,x_4)$ satisfies $Q=0$, that is, the six points are on a conic. Since a (nonsingular) conic is isomorphic to a line, such arrangements form the configuration space  
$$X(2,6)=GL(2)\backslash\{Mat(2,6)\mid \mbox{any }2\times2\mbox{ minor}\not=0\}/(\C^\times)^6$$ of six points on the projective line: if we 
represent a point of $X(2,6)$ by a matrix of the form
\[ z=\left( 
\begin{array}{cccccc}
1&0&1&1& 1& 1\\ 0&1&1&z_1&z_2&z_3
\end{array}\right), \] 
where 
$$\prod_{i=1}^3z_i(z_i-1)\prod_{1\le i<j\le3}(z_i-z_j)\not=0,$$
then the degenerate arrangements in question  can be parametrized by $z=(z_1,z_2,z_3)$ as
$$x_1=(1-z_1)/(1-z_2),\quad x_2=(1-z_1)/(1-z_3),\quad x_3=z_1/z_2,\quad x_4=z_1/z_3.$$
Among the two $S_6$-equivariant projective embedding of $X(2,6)$ presented in 
\cite{Yo1}, let us recall the following one given by the fifteen polynomials
$$D_{ij}(z)D_{kl}(z)D_{mn}(z),\quad\{i,j,k,l,m\}=\{1,\dots,6\},$$
where $D_{ij}(z)$ is the $(i,j)$-minor of the $2\times6$-matrix $z$. Their actual forms are given by
$$(z_i-1)(z_j-z_k), \quad z_j-z_k,\quad z_i(z_j-z_k),\quad z_i(z_j-1).$$
It is known and easy to show that the image is projectively equivalent
to the so-called Segre cubic defined by
$$t_0+\cdots,t_5=0,\quad (t_0)^3+\cdots+(t_5)^3=0.$$

On the other hand, let us prolong the domain of definition of the map $\varphi$
on these degenerate arrangements by the same forty polynomials. 
Then the map $\varphi$ in $z$-coodinates is given by $y_{1}=\cdots y_{10}=0$ and
$$\begin{array}{lll}
cy_{11}=-z_1(z_2-z_3),& cy_{12}=-z_2+z_1,& cy_{13}=z_1-z_3,\\
cy_{14}=-(-1+z_1)(z_2-z_3),& cy_{15}=(-1+z_1)z_3,& cy_{16}=(-1+z_1)z_2,\\
cy_{17}=z_1(-1+z_3),& cy_{18}=z_1(-1+z_2),& cy_{19}=-(-z_2+z_1)z_3,\\
cy_{20}=-(z_1-z_3)z_2,& cy_{21}=-(-z_2+z_1)(-1+z_3),&cy_{22}=-(z_1-z_3)(-1+z_2),\\ 
cy_{23}=-(-1+z_1)(z_2-z_3), &cy_{24}=z_1-z_3,& cy_{25}=-z_2+z_1,\\
cy_{26}=-z_1(z_2-z_3),& cy_{27}=-(z_1-z_3)(-1+z_2),& cy_{28}=-(-z_2+z_1)(-1+z_3),\\
cy_{29}=-(z_1-z_3)z_2,& cy_{30}=-(-z_2+z_1)z_3,& cy_{31}=(-1+z_3)z_2,\\
cy_{32}=(-1+z_2)z_3,& cy_{33}=(-1+z_2)z_3,& cy_{34}=(-1+z_3)z_2,\\
cy_{35}=z_1(-1+z_2),& cy_{36}=z_1(-1+z_3),& cy_{37}=(-1+z_1)z_2,\\
cy_{38}=(-1+z_1)z_3,& cy_{39}=z_2-z_3,& cy_{40}=z_2-z_3,\end{array} $$
where
$$c={(z_1-1)(z_1-z_3)(z_2-z_3)(z_1-z_2)z_1\over(1-z_2)^2z_3^2(1-z_3)^2z_2^2}.$$
This shows that the prolonged $\varphi$ gives exactly the embedding of the arranged arragements given by $\{x\mid Q(x)=0\}$, isomorphic to $X(2,6)$, given above.
\par\medskip
Further put
$$x_1=t\xi_1,\quad x_2=t\xi_2,\quad x_3=t\xi_3,\quad x_4=1+t\xi_4,$$
and let $t$ tends to zero in $\varphi(x)$. We can easily see that the limit is a point whchi is independ of $(\xi_1,\dots,\xi_4)$.
\par\medskip
These results together with the facts on Naruki's cross-ratio $\mathcal C$ 
variety reviewd in \S1.3, we can readily show 
\begin{thm}The closure of the
image of $M$ under $\varphi$ is isomorphic to the variety obtained from Naruki's cross-ratio $\mathcal C$ by
blowing down the 40 exceptional divisors of the second kind to points. \end{thm}
\noindent
{\bf Remark:} This variety is isomorphic to the Stake compactification of the modular variety obtained in \cite{ACT}. Note that the 40 cusps correspond to the 40 points obtained by blowing down the divisors of the second kind.
\section{Defining equations}
We will find a set of generators of the ideal defining the closure $\overline{\varphi(M)}$ of the image of $M$ in $\PP^{40-1}.$ 
Among the forty polynomials $y_1(x),\dots,y_{40}(x)$, we can find some relations. The Pr\"ucker relations (\cite{Yo3}) of the $3\times3$-minors of a $3\times6$-matrix yield linear relations among $y_{1}(x),\dots,y_{10}(x)$; for example,
$$(124,356)-(145,236)+(146,235)-(134,256)=y_3(x)-y_4(x)+y_5(x)-y_6(x)=0.$$
On the other hand, it is easy to check the cubic relation
\begin{eqnarray*}&&(124,356)(15,23,46)(13,26,45)-(145,236)(13,56,24)(12,35,46)
\\&&=y_{3}(x)y_{13}(x)y_{21}(x) - y_{4}(x)y_{14}(x)y_{20}(x) =0.\end{eqnarray*}
Let $V$ be the subvariety of $\PP^{40-1}$, coordinatized by $y_1,\dots,y_{40}$,
defined  by the $G$-orbits of the linear and the cubic equations:
$$y_3-y_4+y_5-y_6=0,\quad y_{3}y_{13}y_{21} - y_{4}y_{14}y_{20}=0.$$
\begin{thm}$\overline{\varphi(M)}=V.$\end{thm}
\subsection{An outline of the proof}
By operating the group $G$ to the linear equation, we get many linear relations among $y_1,\dots,y_{40}$. 
These relations form a system of rank 30, that is, all $y$'s can be expressed linearly in terms of ten chosen ones. For example, in terms of 
$$y_{1},\ y_{3},\ y_{4},\ y_{5},\ y_{7},\ y_{11},\ y_{12},\ y_{13},\ y_{15},\ y_{19},$$
the remaining thirty ones are expressed as 
\begin{eqnarray*}
y_{2} &=& y_{1} - y_{5} + y_{4},\quad y_{6} = -y_{4} + y_{5} + y_{3},\quad y_{8} = -y_{3} - y_{1} + y_{7},\\
y_{9} &=& -y_{1} + y_{7} - y_{4},\quad y_{10} = y_{7} - y_{5} - y_{3},\quad y_{14} = y_{11} + y_{13} - y_{12},\\ 
y_{16} &=& -y_{1} + y_{12} - y_{13} - y_{11} + y_{15},\quad y_{17} = y_{15} - y_{1} - y_{13},\\
y_{18}&=&-y_{13} - y_{11} + y_{15},\quad y_{20} = y_{19} + y_{3} + y_{11},\quad y_{21} = y_{4} + y_{19} + y_{12},\\
y_{22} &=& y_{19} + y_{3} + y_{11} + y_{13} + y_{5},\quad y_{23} = y_{1} - y_{12} + y_{13} + y_{3} + y_{11},\\ 
y_{24} &=& y_{4} + y_{1} + y_{13},\quad y_{25} = y_{12} + y_{5} - y_{1},\\ 
y_{26}&=&-y_{4} + y_{5} - y_{1} + y_{3} + y_{11},\quad y_{27}=y_{19} + y_{4} + y_{1} + y_{3} - y_{7} + y_{11} + y_{13},\\ y_{28} &=& -y_{7} + y_{19} + y_{12} + y_{5} + y_{3},\quad
y_{29} = y_{3} - y_{7} + y_{11} + y_{19} + y_{5},\\
y_{30} &=& y_{3} + y_{1} - y_{7} + y_{19} + y_{4},\quad 
y_{31}= -y_{1} - y_{13} + y_{19} + y_{15} + y_{12},\\
y_{32} &=& y_{19} + y_{3} + y_{15},\quad 
y_{33} = y_{19} + y_{4} + y_{15},\\
y_{34} &=& y_{19} - y_{13} + y_{15} + y_{12} + y_{5} - y_{1} + y_{3},\\
y_{35} &=& -y_{5} - y_{1} - y_{3} + y_{7} - y_{11} + y_{15} - y_{13},\\
y_{36} &=& -y_{1} + y_{7} - y_{4} - y_{13} + y_{15},\\
y_{37} &=& -y_{3} + y_{7} - y_{11} + y_{15} - y_{1} + y_{12} - y_{13},\\ 
y_{38} &=& -y_{1} + y_{7} + y_{15},\quad
y_{39} = y_{1} - y_{12} + y_{13},\\ 
y_{40} &=& -y_{5} + y_{4} + y_{1} + y_{13} - y_{12}.
\end{eqnarray*}

By operating the group $G$ to the cubic equation, we get many such relations.
Substituting the expressions of the $y$'s obtained above, 
we get cubic relations in terms of the chosen ten $y$'s; 
let us here rename the ten coordinates as:
\begin{eqnarray*}&&g_1=y_{1},\ g_2=y_{3},\ g_3=y_{4},\ g_4=y_{5},\ g_5=y_{7},\\
&& g_6=y_{11},\ g_7=y_{12},\ g_8=y_{13},\ g_9=y_{15},\ g_0=y_{19}.\end{eqnarray*}
Among these cubic equations we can find exactly thirty linearly independent ones. 
Therefore the variety $V$ is isomorphic to a subvariety of $\PP^9$, coodinatized by $g_1,\dots,g_9,g_0$, defined by thirty cubic equations, say $cub_1,\dots, cub_{30}.$  Some of them which are relatively simple are shown below:
\begin{eqnarray*}
cub_{1}&:=&   g_2g_8g_0+g_2g_8g_7-g_3g_6g_0-g_2g_3g_6-g_3g_6^2-g_8g_0g_3-g_3g_6g_8\\
&&+g_3g_7g_0+g_2g_3g_7+g_6g_7g_3,\\
cub_{2}&:=&   g_0^2g_1+g_3g_1g_0+g_1^2g_0+g_2g_0g_1-g_5g_1g_0+g_0g_1g_6+g_8
g_1g_0+g_2g_3g_6\\
&&+g_2g_6g_1+g_2g_8g_6+g_3g_1g_6+g_1^2g_6+g_8g_1g_6-g_5g_3g_6
-g_5g_1g_6-g_5g_8g_6,\\
cub_{3}&:=&   -g_8g_0g_3+g_8g_0g_4+g_2g_8g_0-g_3g_6g_0-g_2g_3g_6-g_3g_6^2-g_3g_6g_8-g_3g_6g_4,\\
cub_{4}&:=&   g_2g_0g_7+g_2^2g_7+g_2g_6g_7+g_2g_8g_7+g_2g_7g_4-g_0g_4g_6
-g_8g_0g_4+g_0g_4g_7,\\
cub_{5}&:=&   -g_5g_3g_2-g_2g_5g_1-g_5g_8g_2+g_2g_3g_0+g_2g_0g_1+g_2g_8g_0+
g_2g_3g_7\\
&&+g_2g_1g_7+g_2g_8g_7+g_2g_3g_4+g_2g_1g_4+g_2g_8g_4+g_2^2g_3+g_2^2g_1+g_2^2g_8+g_5g_4g_0\\
&&-g_5g_3g_0-g_5g_1g_0-g_5g_8g_0+g_5g_7g_0,\\
cub_{6}&:=&   g_5g_8g_1-g_8g_1g_4-g_2g_8g_1-g_5g_8g_7+g_4g_7g_8+g_2g_8
g_7\\
&&+g_5g_8^2-g_8^2g_4-g_2g_8^2+g_5g_8g_2-g_2g_8g_4-g_2^2g_8+g_5g_8g_6-g_8g_4
g_6-g_2g_8g_6\\
&&-g_5g_7g_6-g_5g_4g_6+g_5g_1g_6,\\
cub_{7}&:=&   -g_9g_2g_3+g_3g_9g_5-g_3g_9g_6-g_9g_3g_0-g_9g_3g_4+2g_9g_2
g_4-g_4g_9g_5\\
&&+g_4g_9g_6+g_9g_0g_4+g_9g_4^2+g_9g_2^2-g_9g_5g_2+g_9g_2g_6+g_9
g_0g_2\\
&&-g_5g_6g_0-g_5g_3g_6-g_5g_6g_9,\\
cub_{8}&:=&   -g_5g_3g_7+g_5g_4g_7+g_2g_5g_7+g_5g_7g_6-g_2g_3g_1-g_2g_0g_1-
g_2g_1g_7+g_5g_3g_1\\
&&+g_5g_1g_0-g_3g_1g_6-g_0g_1g_6-g_1g_7g_6-g_3g_1g_0-g_0^2
g_1-g_1g_7g_0\\
&&-g_1g_4g_3-g_1g_4g_0-g_1g_4g_7,\\
cub_{9}&:=&   -g_8g_1g_3-g_3g_1g_6+g_3g_8g_7+g_6g_7g_3-g_9g_3g_7-g_3g_8
^2-2g_3g_6g_8\\
&&+g_9g_3g_8-g_3g_2g_8-g_2g_3g_6+g_9g_2g_3-g_3g_6^2+g_3g_9g_6-
g_3g_1g_0-g_2g_3g_1\\
&&-g_0^2g_1-g_2g_0g_1-g_9g_0g_1-g_1g_7g_0-g_2g_1g_7-g_9g_1
g_7,\\
cub_{10}&:=&   g_5g_8g_0-g_8g_0g_4-g_2g_8g_0+g_5g_8g_9-g_9g_8g_4-g_2g_9
g_8+g_3g_6g_8\\
&&+g_8g_0g_3+g_2g_3g_6-g_5g_3g_6+g_3g_6^2+g_3g_6g_0+g_3g_6g_4-
g_9g_2g_3+g_3g_9g_5\\
&&-g_3g_9g_6-g_9g_3g_0-g_9g_3g_4,\\
cub_{11}&:=&   g_2g_5g_1+g_5g_1^2-g_5^2g_1+g_5g_1g_0+g_5g_3g_1-g_2g_5g_7-g_5g_1
g_7+g_5^2g_7\\
&&-g_5g_7g_0-g_5g_3g_7+g_5g_8g_2+g_5g_8g_1-g_5^2g_8+g_5g_8g_0+g_5
g_3g_8\\
&&-g_3g_2g_8-g_2g_8g_0-g_2g_8g_7,\\
&\vdots&\\
cub_{19}&:=&   -g_5g_1g_0-g_2g_5g_1+g_0^2g_1+2g_2g_0g_1+g_0g_1g_6+g_1g_7
g_0+g_2g_1g_7+g_1g_7g_6\\
&&+g_1g_4g_0+g_2g_1g_4+g_1g_4g_6+g_2^2g_1+g_2g_6g_1-g_5g_7
g_6-g_5g_4g_6,\\
&\vdots&
\end{eqnarray*}

We first study this system over the field $K:=\C(g_1,\dots, g5)$. 
Geometrically speaking, we project the variety $V$ onto the 4-dimensional 
projective space coordinatized by $g_1:\cdots:g_5$, and study the generic 
fibre of the projection
$$\pi: \PP^9\supset V\ni g_1:\cdots:g_0\longmapsto g_1:\cdots:g_5\in \PP^4.$$
We shall prove that $\pi$ is generically two-to-one. 
This implies $\pi$ is two-to-one on
$$V^\circ:=V-V\cap\cup_{j=1}^{40}\{y_j=0\}.$$
Thus the argument in \S2.4 shows that $\varphi:M\to V^\circ$ is an isomorphism.

\par
We next study the intersection $V\cap\{y_1=0\}$, and prove that $V\cap\cup_{j=1}^{40}\{y_j=0\}$ is the totality of the $G$-orbit of the closure of the image of 
$X(2,6)$ under the prolonged $\varphi$, stated in \S2.5.
\subsection{Computation over $K$}
 From $cub_{3}=0$ and $cub_{7}=0,$ we can solve $g_8$ and $g_9$ as:
\begin{eqnarray*}
g_8 &=& g_3g_6(g_2+g_6+g_0+g_4)/(g_0g_2-g_0g_3+g_0g_4-g_3g_6),\\
g_9 &=& g_5g_6(g_0+g_3)/(-g_2g_3+g_5g_3-g_3g_6-g_0g_3-g_4g_3+2g_2g_4-g_5g_4+g_4g_6+g_0g_4\\
&&+g_4^2+g_2^2-g_5g_2+g_2g_6+g_0g_2-g_5g_6).
\end{eqnarray*}
Substituting these into $cub_{19}$, we can solve $g_6$:
$$g_6=-{g_1(-g_5g_0-g_5g_2+g_0^2+2g_0g_2+g_0g_4+g_7g_0+g_2g_7+g_2^2+g_2g_4)\over g_0g_1+g_1g_4+g_1g_7+g_2g_1-g_5g_7-g_5g_4}.$$
Substituting these expressions into $cub_{1},\dots,cub_{30}$, we have $cub_{3}=cub_{7}=cub_{19}=0$ of course and $cub_{10}=0$; though most of the remaining ones are complicated, $cub_{8}$ is relatively simple:
$$cub_{8}=-(g_1-g_5)qq/(g_0g_1+g_1g_4+g_1g_7+g_2g_1-g_5g_7-g_5g_4),$$
where 
\begin{eqnarray*}
qq&=&g_1g_4t^2+(-g_2g_5-g_5g_4+g_1g_4+g_5g_3)s^2+2g_1g_4st\\
&+&(-g_5g_3g_1+g_2g_1g_5+g_1g_4g_3-g_4^2g_5+g_1g_2g_4+g_5g_3g_4-g_2g_5g_4+g_4^2g_1)s\\
&+&(-g_4g_5g_1+g_4^2g_1+g_1g_4g_3+g_1g_2g_4)t
+g_3g_2g_1g_4-g_5g_3g_1g_4+g_1g_4^2g_3
\end{eqnarray*}
is a quadratic form in 
$$t:=g_0,\quad s:=g_7.$$
By this quadratic equation, we reduce the degree with respect to $t$ of the 
numerators of the $cub$'s to 1. In this way we get 26 $(=30-4)$ equations of the form
$$cube_j:a_jt+b_j=0,\quad j=1,\dots,30,\quad j \not=3,7,10,19,$$
where $a_j,b_j\in K[s].$ The polynomials
$$D_{jk}=a_jb_k-a_kb_j\in K[s],\quad j<k$$
have a unique common factor $dd$, which is quadratic in $s$; its actual form is given by
\begin{eqnarray*}dd&:=&
(g_1^2g_4^2+2g_2g_1^2g_4+g_3^2g_5^2+g_2^2g_3^2+g_2^2g_1^2+g_5^2g_4^2+2g_3g_2g_1g_4-2g_2g_1g_5g_4\\
&&-2g_3g_5^2g_4-2g_1g_5g_4^2+2g_2^2g_3g_1-2g_3^2g_5g_2+2g_5g_3g_2g_4
+2g_5g_3g_1g_4-2g_5g_3g_1g_2)s^2\\
&+&
(-g_2^2g_1^3-2g_1^3g_2g_4+2g_1^2g_2g_4^2+2g_1^2g_5g_4^2+g_1^2g_4^3-g_1^3g_4^2-2g_1g_5g_4^3\\
&&-2g_3g_5^2g_4^2-2g_3g_2^2g_1^2+2g_3g_2g_1g_4^2+g_3^2g_2^2g_4+g_5^2g_4^3+
2g_1^2g_2g_5g_4\\
&&+2g_2^2g_3g_1g_4-4g_3g_5g_1g_2g_4+g_5^2g_3^2g_4+g_2^2g_1^2g_4-g_1
g_5^2g_4^2-g_2^2g_3^2g_1\\
&&-g_3^2g_5^2g_1+2g_5g_3g_2g_4^2+2g_5g_1^2g_2g_3+2g_3^2g_5
g_1g_2-2g_1g_2g_5g_4^2\\
&&+2g_3g_5^2g_1g_4+2g_3g_1g_5g_4^2-2g_3g_5g_1^2g_4-2g_3^2
g_5g_2g_4-2g_1^2g_2g_3g_4)s\\
&-&g_1^2g_3^2g_4^2+g_3g_1^2g_4^3-g_1^3g_3g_4^2-g_1^2g_2g_3^2g_4-g_1^3g_3g_2g_4-g_3^2g_1g_2^2g_4\\
&&-g_3^2g_1g_2g_4^2+g_3g_1g_2g_4^3+g_2^2g_3g_1g_4^2-g_3g_5^2g_1^2g_4+g_1^3
g_3g_5g_4\\
&&+g_3^2g_5g_1g_4^2-g_3^2g_5^2g_1g_4+g_3g_5^2g_1g_4^2-g_3g_5g_1g_4^3-g_3g_1^2
g_2^2g_4\\
&&+g_3^2g_5g_1^2g_4-2g_3g_1g_2g_5g_4^2+2g_3g_1^2g_2g_5g_4+2g_3^2g_1g_2g_5
g_4.\end{eqnarray*}
So $dd(s)=0$ is the compatibility condition of the system of equations $cube_j:a_j(s)t+b_j(s)=0$, among which 
$cube_{11}:a_{11}(s)t+b_{11}(s)=0$
has the minimal degree of coefficients; actually, degrees of $a_{11}$ and $b_{11}$ with respect to $s$ are 3 and 4.

Now we get three equations $dd=0, qq=0$ and $cube_{11}=0$. Substituting the expression $t=-b_{11}/a_{11}$ into $qq$, we get a polynomial in $s$. We can prove that this polynomial has $dd$ as a factor. 
\par\medskip
We have proved that for a given generic point $g_0:\cdots:g_5\in\PP^4$, the inverse under $\pi$ consists of two points. 
They are obtained as follows:  solve the quadratic equation $dd=0$ with respect to $s(=g_7)$. For each root, $t(=g_0)$ is uniquely determined by the linear equation $cube_{11}$. And $g_6,g_8,g_9$ are uniquely determined by $cub3,cub7,cub19$ as stated above. Then all the relations are satisfied.
\subsection{Intersection of $V$ and $\{y_1=0\}$}
We should better work on $V$ living in $\PP^{40-1}$. In terms of the forty coordinates $y_1,\dots,y_{40}$, the cubic equations are 2-term equations. Thus the vanishing of $y_1$ implies that of some other coordinates. Thanks to $G$-action, we can assume that $y_{10}=0$. The vanishing of these two coordinates forces the
vanishing of other coordinates. Tedious case-by-case study shows that every component of $V\cap\{y_1=0\}$ is included in the $G$-orbit of 
$$V\cap\{y_1=\cdots=y_9=y_0=0\},$$
which is the closure of the image of $\{x\mid Q(x)=0\}$ under the prolonged $\varphi$.

\end{document}